\documentclass[11pt]{article}

\usepackage[margin=2.7cm]{geometry}
\usepackage{amsmath,amssymb,amsthm}
\usepackage[round]{natbib}
\usepackage{xcolor}
\usepackage[colorlinks=true,linkcolor=blue!60!black,citecolor=blue!60!black,urlcolor=blue!60!black]{hyperref}

\newtheorem{theorem}{Theorem}
\newtheorem{proposition}{Proposition}
\newtheorem{lemma}{Lemma}
\theoremstyle{definition}
\newtheorem{definition}{Definition}
\newtheorem{assumption}{Assumption}
\theoremstyle{remark}
\newtheorem{remark}{Remark}

\newcommand{\RR}{\mathbb{R}}
\newcommand{\NN}{\mathbb{N}}
\newcommand{\HH}{\mathbb{H}}
\newcommand{\FF}{\mathbb{F}}
\newcommand{\EE}{\mathbb{E}}
\newcommand{\dd}{\mathrm{d}}
\newcommand{\T}{\top}
\newcommand{\bbeta}{\boldsymbol{\beta}}
\newcommand{\vta}{\vartheta}
\newcommand{\eps}{\varepsilon}
\newcommand{\Stief}{\mathcal{V}}
\newcommand{\supp}{\operatorname{supp}}
\newcommand{\fa}{\mathfrak{a}}
\newcommand{\fg}{\mathfrak{g}}
\newcommand{\E}[1]{\EE^{(n)}_{#1}}
\newcommand{\risk}[2]{\mathcal{R}^{(n)}_{r,t}\big(#1,#2\big)}
\newcommand{\fl}[1]{\lfloor #1\rfloor}

\title{Pointwise minimax lower bounds under a multi-index constraint:\\ the price of an unknown orientation}
\author{Nora Serdyukova\\[1mm]
\small Departamento de Estad\'{\i}stica, Universidad de Concepci\'on, Chile\\
\small \texttt{nserdyukova@udec.cl}}
\date{}

\begin{document}
\maketitle

\begin{abstract}
We consider nonparametric regression with random design and non-Gaussian noise, where the regression function has the multi-index form $F(x)=f(\Theta x)$ with an unknown link $f:\RR^m\to\RR$ in an anisotropic H\"older class and an unknown matrix $\Theta$ with $m\le d$ orthonormal rows. For estimation of $F$ at a fixed point we prove that the minimax risk is bounded from below by a multiple of $L^{1/(2\gamma+1)}(n^{-1}\ln n)^{\gamma/(2\gamma+1)}$, where $\gamma$ is the harmonic mean of the smoothness indices of the link divided by $m$. For $m<d$ this extends a Gaussian white noise result to regression. The main new case is $m=d$, where no dimension reduction takes place: the logarithmic factor still appears as soon as the link is anisotropic and the orientation of its axes is unknown. In contrast, when the link is isotropic with smoothness index at most one, the class is rotation invariant up to a constant and the classical rate without logarithm is attained. The logarithmic loss is therefore caused by the unknown orientation of the directions in which the function is least smooth.


\par
\end{abstract}

%% =====================================================================
\section{Introduction}\label{sec:intro}
%% =====================================================================

Minimax lower bounds serve as a benchmark for the best accuracy attainable by an arbitrary estimator. We consider estimation of the value $F(t)$ of a regression function $F:\RR^d\to\RR$ at a fixed point $t\in\RR^d$ from independent observations
\begin{equation}\label{eq:model}
    Y_i = F(X_i)+\eps_i,\qquad i=1,\ldots,n,
\end{equation}
with random design and possibly non-Gaussian noise. The accuracy of an estimator $\tilde F(t)$ is measured by the pointwise risk
\[
    \risk{\tilde F}{F}=\Big(\E{F}\big|\tilde F(t)-F(t)\big|^r\Big)^{1/r},\qquad r\ge 1,
\]
and, for a class $\FF$ of functions, a sequence $\psi_n$ is a lower bound for the minimax risk if
\[
    \liminf_{n\to\infty}\ \psi_n^{-1}\inf_{\tilde F}\sup_{F\in\FF}\risk{\tilde F}{F}>0.
\]

For the classical smoothness classes the answer is well known. Over the isotropic H\"older class $\HH_d(\beta,L)$ the minimax rate is $L^{d/(2\beta+d)}n^{-\beta/(2\beta+d)}$ \citep{Stone80,Tsybakov09}, and over the anisotropic class $\HH_d(\bbeta,L)$, $\bbeta=(\beta_1,\ldots,\beta_d)$, it is $L^{1/(2\gamma+1)}n^{-\gamma/(2\gamma+1)}$ with $\gamma^{-1}=\sum_{k=1}^d\beta_k^{-1}$ \citep{KLP01}. In both cases the dimension enters through the effective smoothness ($\beta/d$, respectively $\gamma$), which deteriorates as $d$ grows.

Structural assumptions are the standard remedy. In the multi-index model
\begin{equation}\label{eq:multi-index}
    F(x)=f\big(\theta_1^\T x,\ldots,\theta_m^\T x\big)=f(\Theta x),
\end{equation}
the link $f:\RR^m\to\RR$ and the matrix $\Theta$ with rows $\theta_1,\ldots,\theta_m$ are unknown. This model is classical in dimension reduction \citep{Stone85,Hristache01,DJS08,KLV21}, and it has recently become a benchmark for understanding how neural networks adapt to low-dimensional structure \citep{BauerKohler19,SchmidtHieber20}; see \citet{BrunaHsu25} for a recent survey. Stone's dimensionality reduction principle suggests that the rate should be that of an $m$-dimensional problem.

For estimation at a point, this is true only up to a logarithmic factor. In the single-index case $m=1$, \citet{LS14} showed that the minimax rate is $L^{1/(2\beta+1)}(n^{-1}\ln n)^{\beta/(2\beta+1)}$ and, for $d=2$, constructed an adaptive estimator attaining it. This logarithm is not the usual price for adaptation \citep{Lepskii90,BrownLow96,Gaiffas07}: it appears in the minimax rate itself and is the price for not knowing the index vector. With the integrated $L_2$ risk no such loss occurs \citep{GaiffasLecue07}. Related effects of an unknown structure on pointwise rates were found for multivariate densities under an independence hypothesis \citep{Rebelles15} and for bivariate densities whose coordinates become independent after an unknown rotation \citep{LepskiRebelles20}. For the multi-index model with $m<d$, the lower bound $L^{1/(2\gamma+1)}(n^{-1}\ln n)^{\gamma/(2\gamma+1)}$ was obtained in the Gaussian white noise model by \citet{Serdyukova18}.

The present note makes two contributions. First, Theorem~\ref{th:m<d} establishes the bound for $m<d$ in the regression model \eqref{eq:model} under mild conditions on the noise and the design; for $m=1$ it recovers the lower bound of \citet[Theorem~4]{LS14}. Second, and this is the main new point, we consider $m=d$, where \eqref{eq:multi-index} only says that $f$ is anisotropic with respect to an unknown orthonormal frame. Theorem~\ref{th:m=d} shows that the logarithmic factor appears as soon as the smoothness indices are not all equal, while Proposition~\ref{prop:iso} shows that it disappears for an isotropic link (with $\beta\le 1$). The logarithm is thus caused neither by dimension reduction nor by the pointwise loss alone, but by the unknown orientation of the directions in which the function is least smooth; Remark~\ref{rem:mechanism} explains the mechanism.

Section~\ref{sec:results} states the results, Section~\ref{sec:proofs} contains the proofs, and Section~\ref{sec:discussion} discusses open problems.

%% =====================================================================
\section{Setting and results}\label{sec:results}
%% =====================================================================

\paragraph{Smoothness classes.} For $g:\RR^m\to\RR$ let $D_k^l g$ denote the partial derivative of order $l$ with respect to $z_k$, let $e_1,\ldots,e_m$ be the canonical basis of $\RR^m$, and let $\fl{\beta}$ be the largest integer strictly smaller than $\beta$.

\begin{definition}\label{def:holder}
Let $\bbeta=(\beta_1,\ldots,\beta_m)\in(0,\infty)^m$ and $L>0$. A function $g:\RR^m\to\RR$ belongs to the anisotropic H\"older class $\HH_m(\bbeta,L)$ if, for every $k=1,\ldots,m$, the derivatives $D_k^l g$, $0\le l\le \fl{\beta_k}$, exist and are continuous, and
\begin{gather*}
    \|D_k^l g\|_\infty\le L,\qquad 0\le l\le\fl{\beta_k},\\
    \big|D_k^{\fl{\beta_k}}g(z+\tau e_k)-D_k^{\fl{\beta_k}}g(z)\big|\le L|\tau|^{\beta_k-\fl{\beta_k}}\qquad\forall z\in\RR^m,\ \tau\in\RR.
\end{gather*}
If $\beta_1=\cdots=\beta_m=\beta$ we write $\HH_m(\beta,L)$.
\end{definition}

Let $\Stief_{m,d}=\{V\in\RR^{m\times d}: VV^\T=I_m\}$ be the set of $m\times d$ matrices with orthonormal rows; $\Stief_{d,d}$ is the orthogonal group.

\begin{definition}\label{def:class}
For $1\le m\le d$, the class of multi-index functions is
\[
    \FF_{m,d}(\bbeta,L)=\big\{F:\RR^d\to\RR\ \big|\ F(x)=f(\Theta x),\ \Theta\in\Stief_{m,d},\ f\in\HH_m(\bbeta,L)\big\}.
\]
\end{definition}

For $m=1$ this is the class of single-index functions $F(x)=f(\theta^\T x)$, $\theta\in\mathbb{S}^{d-1}$, studied in \citet{LS13,LS14}. Requiring only linearly independent unit rows in Definition~\ref{def:class} gives a larger class, so the lower bounds below remain valid for it. Throughout,
\[
    \gamma^{-1}=\sum_{k=1}^m\beta_k^{-1},\qquad
    \psi_n(\gamma,L)=L^{1/(2\gamma+1)}\big(n^{-1}\ln n\big)^{\gamma/(2\gamma+1)}.
\]

\paragraph{Assumptions.} In \eqref{eq:model}, $(X_i,\eps_i)$, $i=1,\ldots,n$, are independent copies of $(X,\eps)$, where $X$ and $\eps$ are independent, $\eps$ has a density $p$, and $X$ has a density $g$.

\begin{assumption}\label{ass:noise}
The density $p$ is positive on $\RR$, and there exist $q,Q>0$ such that
\[
    \int_\RR \frac{p(y+\upsilon_1)\,p(y+\upsilon_2)}{p(y)}\,\dd y\le 1+Q|\upsilon_1\upsilon_2|\qquad\forall\,\upsilon_1,\upsilon_2\in[-q,q].
\]
\end{assumption}

For the Gaussian density $\mathcal N(0,\sigma^2)$ the integral equals $\exp(\upsilon_1\upsilon_2/\sigma^2)$, so Assumption~\ref{ass:noise} holds. More generally it holds when $p$ is regular, has finite Fisher information and decreases fast enough at infinity \citep[see the discussion after Assumption~6 in][]{LS14}.

\begin{assumption}\label{ass:design}
There exists $\fg>0$ such that, for every $k\in\{1,\ldots,d\}$ and every $O\in\Stief_{k,d}$, the random vector $OX$ has a density bounded by $\fg$.
\end{assumption}

Taking $k=d$ shows that $\|g\|_\infty\le\fg$. Assumption~\ref{ass:design} holds, for instance, for every bounded density with bounded support, and for every density with $g(x)\le C(1+|x|_2)^{-\varpi}$, $\varpi>d-1$, since the density of $OX$ at $u$ is the integral of $g$ over a $(d-k)$-dimensional affine subspace. Only $k\in\{m,2m,d\}$ is used in the proofs, and only $k=d$, i.e.\ $\|g\|_\infty\le\fg$, is used in Theorem~\ref{th:m=d}.

\paragraph{Results.}

\begin{theorem}[Dimension reduction, $m<d$]\label{th:m<d}
Let Assumptions~\ref{ass:noise} and~\ref{ass:design} hold, $d\ge 2$, $1\le m<d$, $\bbeta\in(0,\infty)^m$, $L>0$, $r\ge1$ and $t\in\RR^d$. There exists $\varkappa>0$ independent of $n$ and $L$ such that, for all $n$ large enough,
\[
    \inf_{\tilde F}\ \sup_{F\in\FF_{m,d}(\bbeta,L)}\risk{\tilde F}{F}\ \ge\ \varkappa\,\psi_n(\gamma,L),
\]
where the infimum is taken over all estimators.
\end{theorem}

\begin{theorem}[No dimension reduction, $m=d$]\label{th:m=d}
Let Assumptions~\ref{ass:noise} and~\ref{ass:design} hold, $d\ge2$, $L>0$, $r\ge 1$, $t\in\RR^d$, and let $\bbeta\in(0,\infty)^d$ have at least two distinct components. Then there exists $\varkappa>0$ independent of $n$ and $L$ such that, for all $n$ large enough,
\[
    \inf_{\tilde F}\ \sup_{F\in\FF_{d,d}(\bbeta,L)}\risk{\tilde F}{F}\ \ge\ \varkappa\,\psi_n(\gamma,L).
\]
\end{theorem}

\begin{proposition}[Isotropic link, $m=d$]\label{prop:iso}
Let $\beta\in(0,1]$ and $L>0$. Then
\[
    \HH_d(\beta,L)\subset\FF_{d,d}(\beta,L)\subset\HH_d\big(\beta,d^{\,1-\beta/2}L\big).
\]
Consequently, for fixed $\beta$ and $L$ the pointwise minimax rates over $\FF_{d,d}(\beta,L)$ and over $\HH_d(\beta,L)$ coincide up to constant factors. In particular, whenever $n^{-\beta/(2\beta+d)}$ is the minimax rate over $\HH_d(\beta,L)$ (for instance, under the standard conditions of \citealp{Stone80}, including a design density bounded away from zero near $t$), it is also the minimax rate over $\FF_{d,d}(\beta,L)$, and no logarithmic factor appears.
\end{proposition}

\begin{remark}
In terms of the dimension, Theorems~\ref{th:m<d} and~\ref{th:m=d} agree with Stone's dimensionality reduction principle: the rate is that of an $m$-dimensional anisotropic problem, the dimension being hidden in $\gamma$, up to the factor $(\ln n)^{\gamma/(2\gamma+1)}$. In the isotropic case $\bbeta=(\beta,\ldots,\beta)$ with $m<d$, Theorem~\ref{th:m<d} reads $\psi_n=L^{m/(2\beta+m)}(n^{-1}\ln n)^{\beta/(2\beta+m)}$; here the logarithm is present for every $\beta$, since a function of $m<d$ variables is never rotation invariant in $\RR^d$.
\end{remark}

\begin{remark}[Mechanism]\label{rem:mechanism}
In the proofs, the hypotheses are bumps of height $\asymp\psi_n$ with widths $h_k$ along the rows of $\Theta$, where $h_k$ is smaller in the directions of lower smoothness, and they differ only by rotations through angles of order $N^{-1}$. Two bumps with widths $h_1<h_2$ in a plane, rotated by an angle $\phi$, overlap roughly in proportion $h_1/(h_2\phi)$. When $h_1/h_2$ is polynomially small in $n$, a polynomial number $N=n^b$ of rotated bumps is almost orthogonal, while the Kullback--Leibler divergence between each of them and $F_0\equiv0$ may be as large as a multiple of $\ln N\asymp\ln n$; this produces the logarithm. If $h_1=h_2$, rotations do not separate the hypotheses, which is why Proposition~\ref{prop:iso} holds. For $m<d$ the role of $h_2$ is played by the infinite width of the ridge function in the directions orthogonal to the rows of $\Theta$.
\end{remark}

\begin{remark}\label{rem:iso-beta}
Proposition~\ref{prop:iso} is stated for $\beta\le 1$, where it follows from an elementary argument. For $\beta>1$ the class $\HH_d(\beta,L)$ in Definition~\ref{def:holder} controls only the pure partial derivatives, and the analogous inclusion requires embedding theorems for anisotropic H\"older--Nikol'skii spaces \citep{Nikolskii75}; we do not pursue this here.
\end{remark}

%% =====================================================================
\section{Proofs}\label{sec:proofs}
%% =====================================================================

\subsection{A general lower bound}

For functions $U,G$ set $Z(U,G)=\prod_{i=1}^n p\big(Y_i-U(X_i)\big)/p\big(Y_i-G(X_i)\big)$, the likelihood ratio in model \eqref{eq:model}; the design density cancels. The following lemma was established by \citet{KLP08}; we use it in the form stated in \citet[Lemma~4]{LS14}.

\begin{lemma}\label{lem:KLP}
Let $\FF$ be a nonempty class of functions and $t\in\RR^d$. Assume that, for all $n$ large enough, there exist $N_n\in\NN$, $c>1$, $\lambda_n>0$ and $F_0,F_1,\ldots,F_{N_n}\in\FF$ such that
\begin{gather}
    |F_j(t)-F_0(t)|=\lambda_n,\qquad j=1,\ldots,N_n, \label{eq:lem1}\\
    \E{F_0}\Big(\frac{1}{N_n}\sum_{j=1}^{N_n}Z(F_j,F_0)\Big)^2\le c. \label{eq:lem2}
\end{gather}
Then, for every $r\ge1$ and all $n$ large enough,
\[
    \inf_{\tilde F}\sup_{F\in\FF}\Big(\E{F}\big|\tilde F(t)-F(t)\big|^r\Big)^{1/r}\ \ge\ \frac12\Big[1-\sqrt{\frac{c-1}{c+3}}\,\Big]\lambda_n.
\]
\end{lemma}

\subsection{Common construction}\label{sec:construction}

Fix $m\le d$ and $\bbeta\in(0,\infty)^m$, and set
\begin{equation}\label{eq:exponents}
    a=\frac{2\gamma}{2\gamma+1},\qquad a_k=\frac{a}{\beta_k},\quad k=1,\ldots,m,
    \qquad\text{so that}\qquad a=\beta_k a_k,\quad 2a+\sum_{k=1}^m a_k=2.
\end{equation}
Indeed, $\sum_k a_k=a/\gamma$ and $2a+a/\gamma=a(2\gamma+1)/\gamma=2$.

\paragraph{The bump.} Let $w\in C^\infty(\RR)$ with $\supp(w)\subset(-1/2,1/2)$ and $w(0)\neq0$, and set $W(z)=c_w\prod_{k=1}^m w(z_k)$, $z\in\RR^m$. For $l\le\fl{\beta_k}$, $D_k^lW(z)=c_w w^{(l)}(z_k)\prod_{i\ne k}w(z_i)$, and since $|w^{(l)}(u+\tau)-w^{(l)}(u)|\le \min\{2\|w^{(l)}\|_\infty,\|w^{(l+1)}\|_\infty|\tau|\}$, every such derivative is H\"older continuous in $z_k$ with any exponent in $(0,1]$. Hence $c_w>0$ can be chosen, depending only on $w$ and $\bbeta$, such that $W\in\HH_m(\bbeta,1)$.

\paragraph{Scaling.} Let $\fa>0$ be a constant to be chosen (independent of $n$ and $L$) and set
\[
    \delta=\fa L^{-1}\sqrt{n^{-1}\ln n},\qquad A=\delta^{a},\qquad h_k=\delta^{a_k},\quad k=1,\ldots,m.
\]
For $u\in\RR^m$ write $h^{-1}u=(u_1/h_1,\ldots,u_m/h_m)$. For $V\in\Stief_{m,d}$ define
\begin{equation}\label{eq:FV}
    F_V(x)=f_V(Vx),\qquad f_V(z)=L A\,W\big(h^{-1}(z-Vt)\big).
\end{equation}

\noindent\emph{(i) Membership.} Let $n$ be so large that $\delta\le1$. For $l\le\fl{\beta_k}$, $\|D_k^lf_V\|_\infty\le LAh_k^{-l}=L\delta^{a_k(\beta_k-l)}\le L$, and, since $Ah_k^{-\beta_k}=1$ by \eqref{eq:exponents},
\[
    \big|D_k^{\fl{\beta_k}}f_V(z+\tau e_k)-D_k^{\fl{\beta_k}}f_V(z)\big|\le LAh_k^{-\fl{\beta_k}}\Big|\frac{\tau}{h_k}\Big|^{\beta_k-\fl{\beta_k}}=L|\tau|^{\beta_k-\fl{\beta_k}}.
\]
Thus $f_V\in\HH_m(\bbeta,L)$ and $F_V\in\FF_{m,d}(\bbeta,L)$.

\noindent\emph{(ii) Separation.} By \eqref{eq:exponents},
\begin{equation}\label{eq:lambda}
    F_V(t)=LA\,W(0)=W(0)\,\fa^{a}L^{1-a}\big(n^{-1}\ln n\big)^{a/2}=W(0)\,\fa^{a}\,\psi_n(\gamma,L),
\end{equation}
since $1-a=1/(2\gamma+1)$ and $a/2=\gamma/(2\gamma+1)$.

\noindent\emph{(iii) Size.} $\|F_V\|_\infty=LA\|W\|_\infty\to0$, and by \eqref{eq:exponents}
\begin{equation}\label{eq:L2A2h}
    L^2A^2\prod_{k=1}^m h_k=L^2\delta^{\,2a+\sum_k a_k}=L^2\delta^2=\fa^2\,n^{-1}\ln n.
\end{equation}

\paragraph{Reduction to cross terms.} Let $b>0$ be chosen below, $N=N_n=\lfloor n^b\rfloor$, and let $V_1,\ldots,V_N\in\Stief_{m,d}$ be specified in the proofs of the theorems. Put $F_0\equiv0\in\FF_{m,d}(\bbeta,L)$ and $F_j=F_{V_j}$, $j=1,\ldots,N$. By \eqref{eq:lambda}, condition \eqref{eq:lem1} holds with $\lambda_n=|W(0)|\fa^a\psi_n(\gamma,L)$.

Under $F_0$ we have $Y_i=\eps_i$. Using independence and Assumption~\ref{ass:noise}, which applies as soon as $\|F_j\|_\infty\le q$, we obtain for all $i,j$ (including $i=j$)
\begin{equation}\label{eq:EZZ}
    \E{F_0}\big[Z(F_i,F_0)Z(F_j,F_0)\big]=\Big\{\EE\!\int_\RR\frac{p(y-F_i(X))\,p(y-F_j(X))}{p(y)}\dd y\Big\}^n
    \le\exp\big\{nQ\,\EE|F_i(X)F_j(X)|\big\}.
\end{equation}

\noindent\emph{Diagonal terms.} By Assumption~\ref{ass:design} with $k=m$ and $O=V_j$, the density of $V_jX$ is bounded by $\fg$, so by \eqref{eq:L2A2h}
\[
    \EE F_j^2(X)=\EE f_{V_j}^2(V_jX)\le\fg\int_{\RR^m}f_{V_j}^2(z)\,\dd z=\fg\|W\|_2^2L^2A^2\prod_{k=1}^mh_k=\fg\|W\|_2^2\,\fa^2\,n^{-1}\ln n,
\]
where $\|W\|_2$ is the $L_2(\RR^m)$ norm. Hence $\E{F_0}Z^2(F_j,F_0)\le n^{Q\fg\|W\|_2^2\fa^2}$. We choose
\begin{equation}\label{eq:choice-a}
    \fa^2=\frac{b}{2Q\fg\|W\|_2^2},
\end{equation}
so that, since $N\ge n^b/2$ for $n$ large, $N^{-2}\sum_{j}\E{F_0}Z^2(F_j,F_0)\le 2n^{-b/2}\to0$.

\noindent\emph{Conclusion.} Consequently,
\[
    \E{F_0}\Big(\frac1N\sum_{j=1}^NZ(F_j,F_0)\Big)^2\le 2n^{-b/2}+\max_{i\ne j}\exp\big\{nQ\,\EE|F_i(X)F_j(X)|\big\}.
\]
If $b$ is chosen so that
\begin{equation}\label{eq:cross-goal}
    \max_{1\le i\ne j\le N}\ n\,\EE|F_i(X)F_j(X)|\ \longrightarrow\ 0,\qquad n\to\infty,
\end{equation}
then \eqref{eq:lem2} holds with $c=2$ for all $n$ large enough, and Lemma~\ref{lem:KLP} yields both theorems with
\[
    \varkappa=\tfrac12\big(1-5^{-1/2}\big)|W(0)|\,\fa^{a},
\]
which depends on $b$, $Q$, $\fg$, $w$ and $\bbeta$, but not on $n$ or $L$. How large $n$ must be does depend on $L$, through $\delta\le1$, $\|F_j\|_\infty\le q$ and the convergence in \eqref{eq:cross-goal}. It remains to choose $V_1,\ldots,V_N$ and $b$ and to verify \eqref{eq:cross-goal}.

Throughout, $\varphi_j=j/N$, $c_j=\cos\varphi_j$, $s_j=\sin\varphi_j$. For $i\ne j$ we have $1/N\le|\varphi_i-\varphi_j|\le1<\pi/2$, hence
\begin{equation}\label{eq:sin}
    |\sin(\varphi_i-\varphi_j)|\ge\frac{2}{\pi}|\varphi_i-\varphi_j|\ge\frac{2}{\pi N}.
\end{equation}
We also write $\|\cdot\|_1$, $\|\cdot\|_\infty$ for the norms of $w$ on $\RR$, and $y=x-t$.

\subsection{Proof of Theorem~\ref{th:m<d}}

Let $1\le m<d$ and $p=m\wedge(d-m)$. Define $V_j\in\Stief_{m,d}$ by its rows
\[
    \vta_{j,k}=c_je_k+s_je_{m+k},\quad k=1,\ldots,p;\qquad \vta_{j,k}=e_k,\quad k=p+1,\ldots,m,
\]
where now $e_1,\ldots,e_d$ is the canonical basis of $\RR^d$. The rows are orthonormal, since for $k\le p$ the indices $k$ and $m+k\le d$ are distinct and different from those used by the other rows. We have
\[
    |F_i(x)F_j(x)|=L^2A^2c_w^2\prod_{k=1}^m\Big|w\Big(\frac{\vta_{j,k}^\T y}{h_k}\Big)\Big|\ \prod_{k=1}^m\Big|w\Big(\frac{\vta_{i,k}^\T y}{h_k}\Big)\Big|.
\]

\noindent\emph{Case $d/2\le m<d$, so that $p=d-m$.} Bound the factors $|w(\vta_{i,k}^\T y/h_k)|$, $k>p$, by $\|w\|_\infty$, and let $M\in\RR^{d\times d}$ have rows $\vta_{j,1},\ldots,\vta_{j,m},\vta_{i,1},\ldots,\vta_{i,p}$. After a permutation of coordinates, $M$ is block diagonal with the identity on $e_{p+1},\ldots,e_m$ and $p$ blocks $\big(\begin{smallmatrix}c_j&s_j\\ c_i&s_i\end{smallmatrix}\big)$ on the pairs $(e_k,e_{m+k})$, $k\le p$; hence $|\det M|=|\sin(\varphi_i-\varphi_j)|^p$. Since $\|g\|_\infty\le\fg$, the change of variables $v=My$ gives
\begin{align*}
    \EE|F_i(X)F_j(X)|&\le L^2A^2c_w^2\|w\|_\infty^{m-p}\,\frac{\fg}{|\det M|}\int_{\RR^d}\prod_{k=1}^m\Big|w\Big(\frac{v_k}{h_k}\Big)\Big|\prod_{k=1}^p\Big|w\Big(\frac{v_{m+k}}{h_k}\Big)\Big|\dd v\\
    &=\frac{c_w^2\fg\|w\|_\infty^{m-p}\|w\|_1^{m+p}}{|\det M|}\,L^2A^2\prod_{k=1}^mh_k\prod_{k=1}^ph_k.
\end{align*}

\noindent\emph{Case $1\le m<d/2$, so that $p=m$.} All rows of $V_i$ and $V_j$ lie in the span of $e_1,\ldots,e_{2m}$. Hence the $2m\times d$ matrix $U$ with rows $\vta_{j,1},\ldots,\vta_{j,m},\vta_{i,1},\ldots,\vta_{i,m}$ factorizes as $U=BO$, where $O\in\Stief_{2m,d}$ has rows $e_1,\ldots,e_{2m}$ and $B\in\RR^{2m\times2m}$ is, after a permutation, block diagonal with the $m$ blocks above, so $|\det B|=|\sin(\varphi_i-\varphi_j)|^m$. By Assumption~\ref{ass:design} with $k=2m$, the density of $UX=B(OX)$ is bounded by $\fg/|\det B|$, and therefore
\[
    \EE|F_i(X)F_j(X)|\le L^2A^2c_w^2\,\frac{\fg}{|\det B|}\int_{\RR^{2m}}\prod_{k=1}^m\Big|w\Big(\frac{v_k}{h_k}\Big)w\Big(\frac{v_{m+k}}{h_k}\Big)\Big|\dd v
    =\frac{c_w^2\fg\|w\|_1^{2m}}{|\det B|}\,L^2A^2\Big(\prod_{k=1}^mh_k\Big)^2.
\]

\noindent\emph{Both cases.} By \eqref{eq:sin} and \eqref{eq:L2A2h}, with $C_1=c_w^2\fg(\pi/2)^p\|w\|_\infty^{m-p}\|w\|_1^{m+p}$ and $s_p=\sum_{k=1}^pa_k$,
\[
    n\,\EE|F_i(X)F_j(X)|\le C_1N^p\,\fa^2\ln n\prod_{k=1}^ph_k
    \le C_1\fa^2\Big(\frac{\fa}{L}\Big)^{s_p}n^{\,bp-s_p/2}(\ln n)^{1+s_p/2}.
\]
Choosing $0<b<s_p/(2p)$, the right-hand side tends to zero uniformly in $i\ne j$, which is \eqref{eq:cross-goal}. This proves Theorem~\ref{th:m<d}.\qed

\subsection{Proof of Theorem~\ref{th:m=d}}

Relabelling the arguments of $f$, together with the rows of $\Theta$, shows that $\FF_{d,d}(\bbeta,L)$ does not change when the components of $\bbeta$ are permuted. Hence we may assume $\beta_1<\beta_2$, so that $a_1>a_2$ and $h_1<h_2$ for $\delta<1$: the bump is narrower in the less smooth direction.

Define $V_j\in\Stief_{d,d}$ as the rotation by $\varphi_j$ in the plane of $e_1,e_2$:
\[
    \vta_{j,1}=c_je_1+s_je_2,\qquad\vta_{j,2}=-s_je_1+c_je_2,\qquad\vta_{j,k}=e_k,\quad k=3,\ldots,d.
\]
In the product $|F_i(x)F_j(x)|$ we keep the factors of $F_j$ with $k\ne2$ and the factor of $F_i$ with $k=1$, and bound the remaining $d$ factors by $\|w\|_\infty$. Thus the two kept factors with $k=1$ correspond to the nearly parallel directions $\vta_{j,1}$ and $\vta_{i,1}$. The matrix $M$ with rows $\vta_{j,1},\vta_{i,1},e_3,\ldots,e_d$ satisfies $|\det M|=|c_js_i-s_jc_i|=|\sin(\varphi_i-\varphi_j)|$. Using $\|g\|_\infty\le\fg$, the change of variables $v=My$ and \eqref{eq:sin}, we obtain
\[
    \EE|F_i(X)F_j(X)|\le\frac{c_w^2\fg\|w\|_\infty^d\|w\|_1^d}{|\det M|}\,L^2A^2h_1^2\prod_{k=3}^dh_k
    \le C_2N\,L^2A^2\prod_{k=1}^dh_k\cdot\frac{h_1}{h_2},
\]
with $C_2=(\pi/2)c_w^2\fg\|w\|_\infty^d\|w\|_1^d$. By \eqref{eq:L2A2h} and $h_1/h_2=\delta^{a_1-a_2}$,
\[
    n\,\EE|F_i(X)F_j(X)|\le C_2\fa^2\Big(\frac{\fa}{L}\Big)^{a_1-a_2}n^{\,b-(a_1-a_2)/2}(\ln n)^{1+(a_1-a_2)/2}.
\]
Since $a_1>a_2$, any $0<b<(a_1-a_2)/2$ gives \eqref{eq:cross-goal}. This proves Theorem~\ref{th:m=d}.\qed

If $\beta_1=\cdots=\beta_d$, then $h_1=h_2$ and the last bound is of order $N\ln n$, which does not tend to zero for any $b>0$. This is consistent with Proposition~\ref{prop:iso}.

\subsection{Proof of Proposition~\ref{prop:iso}}

The first inclusion follows by taking $\Theta=I_d$. For the second, let $\beta\le1$, so that $\fl\beta=0$, and let $F(x)=f(\Theta x)$ with $f\in\HH_d(\beta,L)$ and $\Theta$ orthogonal. Clearly $\|F\|_\infty\le L$. For $u\in\RR^d$, changing one coordinate at a time and using the inequality $\sum_{k=1}^d|u_k|^\beta\le d^{1-\beta/2}|u|_2^\beta$ (H\"older's inequality, $\beta\le2$), we get
\[
    |f(z+u)-f(z)|\le L\sum_{k=1}^d|u_k|^\beta\le d^{1-\beta/2}L\,|u|_2^\beta.
\]
Applying this with $z=\Theta x$ and $u=\tau\Theta e_k$, $|u|_2=|\tau|$, gives $|F(x+\tau e_k)-F(x)|\le d^{1-\beta/2}L|\tau|^\beta$, that is, $F\in\HH_d(\beta,d^{1-\beta/2}L)$. The remaining assertions follow because the minimax risk is monotone with respect to inclusion of classes.\qed

%% =====================================================================
\section{Discussion}\label{sec:discussion}
%% =====================================================================

\paragraph{Sharpness.} For $m=1$ and $d=2$ the bound of Theorem~\ref{th:m<d} is attained by an adaptive estimator \citep{LS14}, and in the density model of \citet{LepskiRebelles20} the analogous bound for an unknown rotation in $\RR^2$ is attained as well. For $m\ge2$, and in particular for $m=d$, matching upper bounds are open. A natural candidate is a kernel estimator with anisotropic bandwidths along orthonormal frames from a grid of polynomial cardinality in $n$, combined with a data-driven selection rule in the spirit of \citet{GL09}; the logarithm in Theorems~\ref{th:m<d} and~\ref{th:m=d} is precisely the cost of searching such a grid. The construction of this estimator is the subject of further work.

\paragraph{Loss function.} The hypotheses in the proofs differ only in a neighbourhood of $t$, and the argument relies on the pointwise loss. For integrated losses no logarithmic factor is expected; for the single-index model with $L_2$ loss this was shown by \citet{GaiffasLecue07}. Thus the logarithm is caused by the combination of a pointwise loss with an unknown anisotropic orientation.

\paragraph{Smoother isotropic links.} Proposition~\ref{prop:iso} covers $\beta\le1$. Extending it to $\beta>1$ amounts to showing that the classes of Definition~\ref{def:holder} with equal exponents are rotation invariant up to constants; see Remark~\ref{rem:iso-beta}.

%% =====================================================================

\end{document}